\documentclass[12pt]{article}
\usepackage{graphicx}

\usepackage{amsmath,amsfonts,amsthm,amssymb}
\newcommand{\homo}{{\mathrm H}}

\begin{document}
\title
{Embedding infinite cyclic covers of knot spaces into
3-space\footnote{All authors are partially supported by a MOSTC
grant and a MOEC grant. The second author is partially supported
by the Centennial fellowship of Princeton University. Part of the
paper is revised from the Master thesis of the second named author
submitted to Peking University.}}
\maketitle

\begin{center}
BOJU JIANG, jiangbj@math.pku.edu.cn
\\
{\it LMAM Depart. of Math.,
Peking University, Beijing  100871, China\\
}

\end{center}

\begin{center}
YI NI,  yni@princeton.edu\\
{\it Depart. of Mathematics,
Princeton University,  NJ  08544, U. S. A.\\
}
\end{center}

\begin{center}
SHICHENG WANG, wangsc@math.pku.edu.cn\\
{\it LMAM Depart. of Math., Peking University, Beijing 100871,
China\\
}
\end{center}

\begin{center}
QING ZHOU, qingzhou@sjtu.edu.cn\\
{\it Depart. of Mathematics,
Jiaotong University, Shanghai 200030, China\\
}
\end{center}

\newpage

\begin{abstract}
We say a knot $k$ in the 3-sphere $\mathbb S^3$ has {\it Property
$IE$} if the infinite cyclic cover of the knot exterior embeds
into $\mathbb S^3$. Clearly all fibred knots have Property $IE$.

There are infinitely many non-fibred knots with Property $IE$ and
infinitely many non-fibred knots without property $IE$. Both kinds
of examples are established here for the first time. Indeed we
show that if a genus 1 non-fibred knot has Property $IE$, then its
Alexander polynomial $\Delta_k(t)$ must be either 1 or
$2t^2-5t+2$, and we give two infinite families of non-fibred genus
1 knots with Property $IE$ and having $\Delta_k(t)=1$ and
$2t^2-5t+2$ respectively.

Hence among genus one non-fibred knots, no alternating knot has
Property $IE$, and there is only one knot with Property $IE$ up to
ten crossings.

We also give an obstruction to embedding infinite cyclic covers of
a compact 3-manifold into any compact 3-manifold.

\vskip 0.2 true cm

 {\it Keywords:} Embedding; Non-fiber knots;
Infinite cyclic coverings.

{\it MSC:} 57M10, 57M25, 57N30
\end{abstract}

\vskip 0.2 true cm
\begin{center}{\bf \S 1. Introduction}\end{center}

In this paper all surfaces and 3-manifolds are orientable, and all
surfaces in 3-manifolds are proper, embedded and two-sided.
Suppose $S$ (resp. $P$) is a surface (resp. 3-manifold) in a
3-manifold $M$, we use $M\setminus S$ (resp. $M\setminus P$) to
denote the manifold obtained by cutting $M$ along $S$ (resp.
removing $\text{int} P$, the interior of $P$, from $M$).

 Suppose $S$ is a connected non-separating
surface in $M$. Then $X=M\setminus S$ has two copies of $S$ in
$\partial X$, denoted by  $S^+\sqcup S^-$. Taking countably many
copies of $X$: $\{X_i\}_{i=-\infty}^{+\infty}$, and identifying
$S^+_{i-1}$ with $S^-_i$ for all $i$, we get an infinite cyclic
cover of $M$, denoted by $\widetilde {M}_S$.

Let $k$ be a knot in $\mathbb S^3$, $E(k)$ be the exterior of $k$,
$S$ be a Seifert surface of $k$. Then $E(k)$ has a unique infinite
cyclic cover, simply denoted by $\widetilde E(k)$. If $k$ is a
fibred knot with fiber $S$, then $\widetilde E(k)$ is homeomorphic
to $S \times\mathbb R$ which clearly embeds into $\mathbb S^3$.
The paper will address the following

{\bf Question 1.} {\it Suppose $k$ is a non-fibred knot, when does
$\widetilde E(k)$ embed into $\mathbb S^3$?}

The third named author was introduced to Question 1 during
conversations with Professor Robert D. Edwards in the spring of
1984, and  Edwards attributed Question 1 to Professor J.
Stallings.

It is natural to ask the following more general and flexible

{\bf Question 2.} {\it When does an infinite cyclic cover of a
compact 3-manifold embed into a compact 3-manifold?}

{\bf Definition 1.1}  We say a knot $k$ in $\mathbb S^3$ has {\it
Property $IE$}, if the infinite cyclic cover $\widetilde E(k)$
embeds into $\mathbb S^3$.  We say a knot $k$ in $\mathbb S^3$ has
{\it Property $DIE$}, if $(\widetilde E(k), \tau) \subset (\mathbb
S^3, f)$, that is, the deck transformation $\tau$ of $\widetilde
E(k)$ embeds into a dynamical system $f$ on $\mathbb S^3$. (We say
a dynamical system $g$ on a space $P$ {\it embeds} into a
dynamical system $f$ on a space $Y$, denoted by $(P, g)\subset (Y,
f)$, if there is an embedding $P\subset Y$ such that $f|P=g$.)

The organization of this paper goes as below.

\S 2 and \S 3 are the main parts of the paper. All knots involved
in \S 2 and \S 3 are of genus 1 and non-fibred. It is well known
that the only genus 1 fibred knots are $3_1$ and $4_1$ in the knot
table.

In \S 2, we give a partial positive answer to Questions 1 and 2.
In {\S 2.1}, beginning with a discrete dynamical system $f$ on
$\mathbb S^3$ (or a compact 3-manifold $Y$), we construct a
compact 3-manifold $M$ (closed or with torus boundary) such that
$(\widetilde M_S, \tau)\subset (\mathbb S^3, f)$ or $\subset (Y,
f)$, where $\tau$ is the deck transformation on the infinite
cyclic cover $\widetilde M_S$. In {\S 2.2} we prove that the
simplest non-trivial example provided by construction in {\S 2.1}
is $E(9_{46})$, the exterior of the 46-th knot of nine crossings
in the knot table, see \cite{R} or \cite{BZ}, therefore provide
the first known positive example to Question 1. A subtle point in
the verification is to choose a right projection of $9_{46}$,
which significantly simplifies the process. But a key point is to
choose $9_{46}$ among all knots in $\mathbb S^3$ to compare with.
In {\S 2.3}, we give  a sufficient condition for the 3-manifolds
constructed in {\S 2.1} to be complements of knots in $\mathbb
S^3$, and then we prove that there are infinitely many non-fibred
genus 1 knots having Property $DIE$ by invoking Thurston and
Soma's results on Gromov volume of 3-manifolds.

In \S 3, we give a partial negative answer to Question 1.  By
invoking Freedman-Freedman's version of Kneser-Haken finiteness
theorem  and results of Gabai (and Novikov) on foliation and on
surgery, we prove that if a genus 1 non-fibred knot $k$ has
Property $IE$, then $E(k)$ is constructed as in \S 2.1, and hence
$k$ has Property $DIE$. It follows that the Alexander polynomial
of such knots must be $1$ or $2t^2-5t+2$, and the Alexander
invariant is also restricted. So ``most" genus 1 non-fibred knots
do not have Property $IE$. In particular, among all non-fibred
genus 1 knots, no alternating knots have Property $IE$, and up to
crossing numbers $\le 10$ only $9_{46}$ has Property $IE$. On the
other hand, two infinite families of genus 1 non-fibred  knots
with Property $IE$ constructed in \S 2.3 have $\Delta_{k}(t)=1$
and $\Delta_{k}(t)=2t^2-5t+2$ respectively.

\S 4 is a remark about Property $IE$ on connected sum, which
provides knots of any given genus $g$ (non-prime when $g>1$), some
of them have Property $IE$ and some do not.

\S 5 gives a homological obstruction to embedding infinite cyclic
covers of a compact 3-manifold into any compact 3-manifold
(Theorem 5.1), therefore gives a partial negative answer to
Question 2.

{\bf Comments.}

1. If we replace the term ``unknotted solid torus" by ``unknotted
handlebody of genus $g$ for any $g>1$", constructions in {\S 2.1}
can be used to study Property $DIE$ of knots with higher genera,
although the arguments become more complicated. The knots having
Property $DIE$ provide interesting dynamics in $\mathbb S^3$.

2. Theorem 5.1 as well as the constructions in \S 2.1  still holds
for closed $n$-manifold  and  connected non-separating bicollared
properly embedded $(n-1)$-submanifold $S$ in $M$.

3. For knot $k$ in $S^3$, the homological obstruction in Theorem
5.1 vanishes for $E(k)$ (read Remark 2). We wonder if Question 2
has positive answer when we restrict to $E(k)$ for knots $k$ in
$\mathbb S^3$.

4. Two references \cite{JNW} and \cite{CL} were not cited in our
proofs. But \cite{JNW} suggested us the construction in \S 2.1 and
\cite{CL} inspired us to prove Lemma 3.1.

{\bf Acknowledgements.}  We are grateful to  Dr. Hao Zheng for
drawing the pictures, to Professor William Browder for a helpful
conversation with the second author, to Professor Robert D.
Edwards for bringing the  third author to this simply stated
intuitive question, to Professor David Gabai for a comment on
alternating knots.

\begin{center}{\bf \S 2. Infinitely many genus 1 non-fibred knots have Property DIE}\end{center}

{\bf\S 2.1.  A construction of compact 3-manifolds having infinite
cyclic covers in $\mathbb S^3$ or in a compact 3-manifold}

Step 1. We first consider a rather general case. Let $Y$ be a
closed 3-manifold, and $P\subset Y$ be a submanifold of dimension
three with connected and non-empty $\partial P$. Suppose that
there is a homeomorphism
$$f: Y\to Y \qquad \text {such that} \qquad f(P) \subset \text{int} P.$$

Let $X=P\setminus f(P)$. Then $\partial X=\partial P\cup
\partial f(P)$. Let $M=X/f$ be the closed $3$-manifold obtained
from $X$ by identifying $\partial P$ and $\partial f(P)$ via $f$,
and $S\subset M$ be the image of $\partial P$ and $\partial f(P)$
after identification. Then $S$ is a connected non-separating
surface in $M$. Clearly the infinite cyclic cover $\widetilde M_S$
is identified with $\cup_{k=-\infty}^{+\infty} f^k(X)\subset Y$
and $f|\cup_{k=-\infty}^{+\infty} f^k(X)$ gives the deck
transformation $\tau$. Hence $(\widetilde M_S, \tau)\subset (Y,
f)$.

We say the construction above is {\it non-trivial}, if $X$ is not
homeomorphic to $\partial P\times [0,1]$.

Step 2. Continue from Step 1. Let $Y=\mathbb S^3$ and let $P$ be
an unknotted solid torus $P$ in $\mathbb S^3$, and let $P'$ be a
solid torus in $\text{int} P$, such that $P'$ is still unknotted
in $\mathbb S^3$. Since both $P$ and $P'$ are unknotted in
$\mathbb S^3$, there is a homeomorphism $f: \mathbb S^3\to\mathbb
S^3$ such that $f(P)=P'$. Then $X=P\setminus f(P)$ is an example
of Step~1.

Step 3. Continue from Step 2. Let $\Gamma$ be a proper arc in $X$
with one end in $\partial P$ and the other in $\partial f(P)$. Let
$N(\Gamma)$ be the regular neighborhood of $\Gamma$ in $X$. Up to
isotopy we may assume $f(\partial P\cap N(\Gamma))=\partial
f(P)\cap N(\Gamma)$. Let $X^*=X\setminus N(\Gamma)$. Then $X^*$ is
obtained from $X$ by digging a tunnel from $\partial P$ to
$\partial f(P)$. Let  $M^*=X^*/f$, $S^*=M^*\cap S$, where $M^*$ is
obtained from $M$ by removing a solid torus. Clearly the infinite
cyclic cover $\widetilde M^*_{S^*}$ is identified with
$\cup_{k=-\infty}^{+\infty} f^k(X^*)\subset\mathbb S^3$, and
$f|\cup_{k=-\infty}^{+\infty} f^k(X^*)$ gives the deck
transformation $\tau$. We summarize the discussion above as

{\bf Proposition 2.1} {\it $M^*$ is a compact 3-manifold with
torus boundary, and $(\widetilde M^*_S, \tau)\subset (\mathbb S^3,
f)$. In particular if $M^*$ is homeomorphic to $E(k)$ for a knot
$k\subset S^3$, then $k$ has Property $DIE$.}

{\bf\S 2.2. The knot $9_{46}$ has Property $DIE$}

A simplest non-trivial construction in Proposition 2.1 is
indicated in Figure 1, where $P'$ is a 2-braid in $P$ and the
tunnel is ``unknotted". In this subsection, all notions in Step 3
of \S 2.1 refer to Figure~1.

\begin{center}%
\includegraphics[totalheight=7cm]{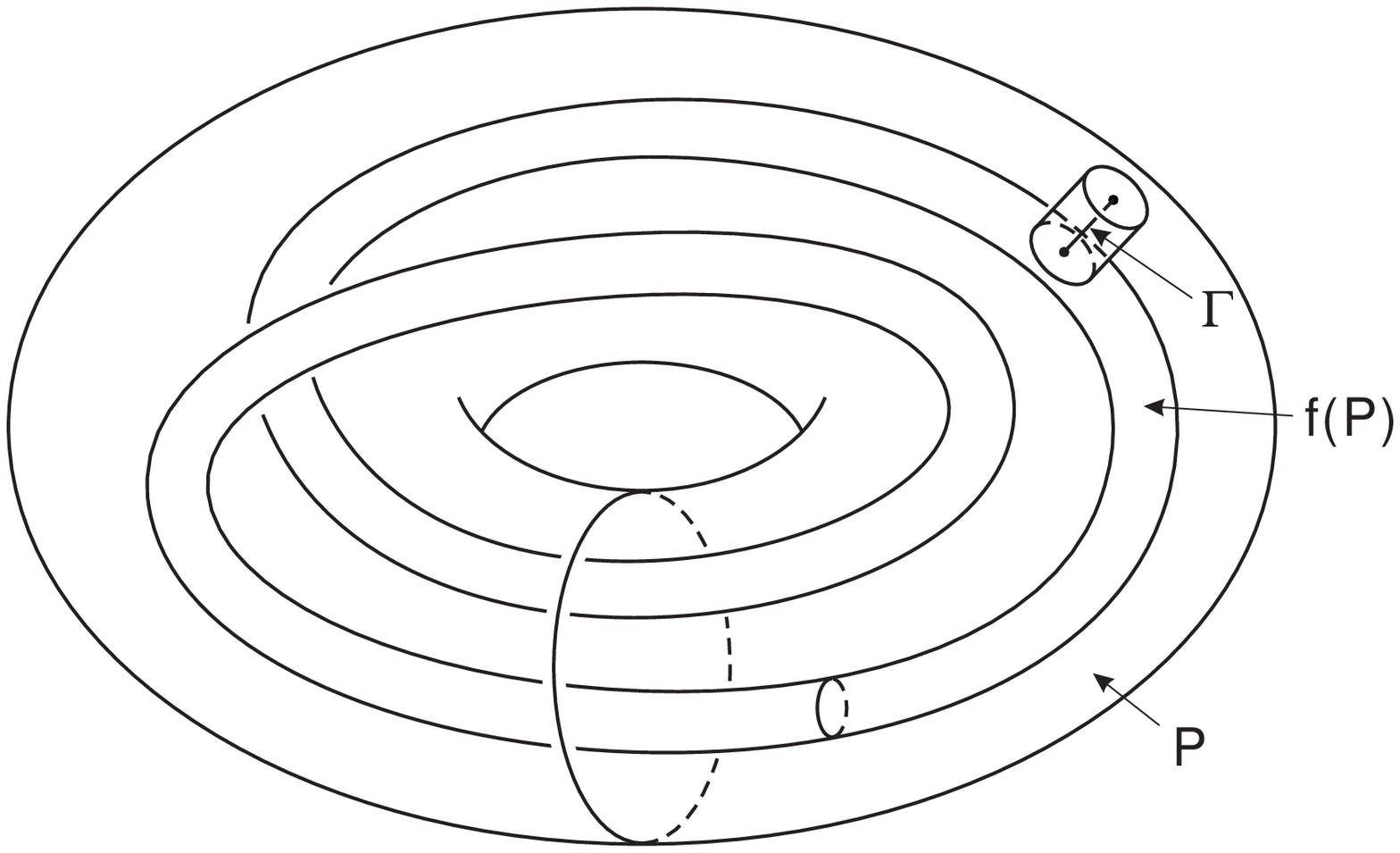}%
\begin{center}%
Figure 1
\end{center}
\end{center}

We will verify  that $M^*$ is homeomorphic to $E=E(9_{46})$, the
exterior of knot $9_{46}\subset\mathbb S^3$ in the knot table. Our
verification consists of three steps:

Step 1. Compute $\pi_1(M^*)$ and $\pi_1(\partial M^*)\subset
\pi_1(M^*)$. Cutting $M^*$ open along $S^*$, we get back to $X^*$,
which is already presented in Figure 1. Its boundary $\partial
X^*=S^*_-\cup {\text{\rm annulus}}\cup S_+^*$, where $S_-^*$ and
$S_+^*$ are 1-punctured tori on the inner boundary $\partial f(P)$
and the outer boundary $\partial P$ respectively. The annulus is
the boundary of the tunnel.

Choose  meridian $\mu_+$ and longitude $\lambda_+$ on $S^*_+$ such
that $\mu_+$ bounds a disc in $P$ and $\lambda_+$ bounds a disc in
$\mathbb S^3\setminus P$. Similarly choose meridian $\mu_-$ and
longitude $\lambda_-$ on $S^*_-$ such that $\mu_-$ bounds a disc
in $f(P)$ and $\lambda_-$ bounds a disc in $\mathbb S^3\setminus
f(P)$, where $\mu_\pm$ and $\lambda_\pm$ are as indicated in
Figure 2.

\begin{center}%
\includegraphics[totalheight=7cm]{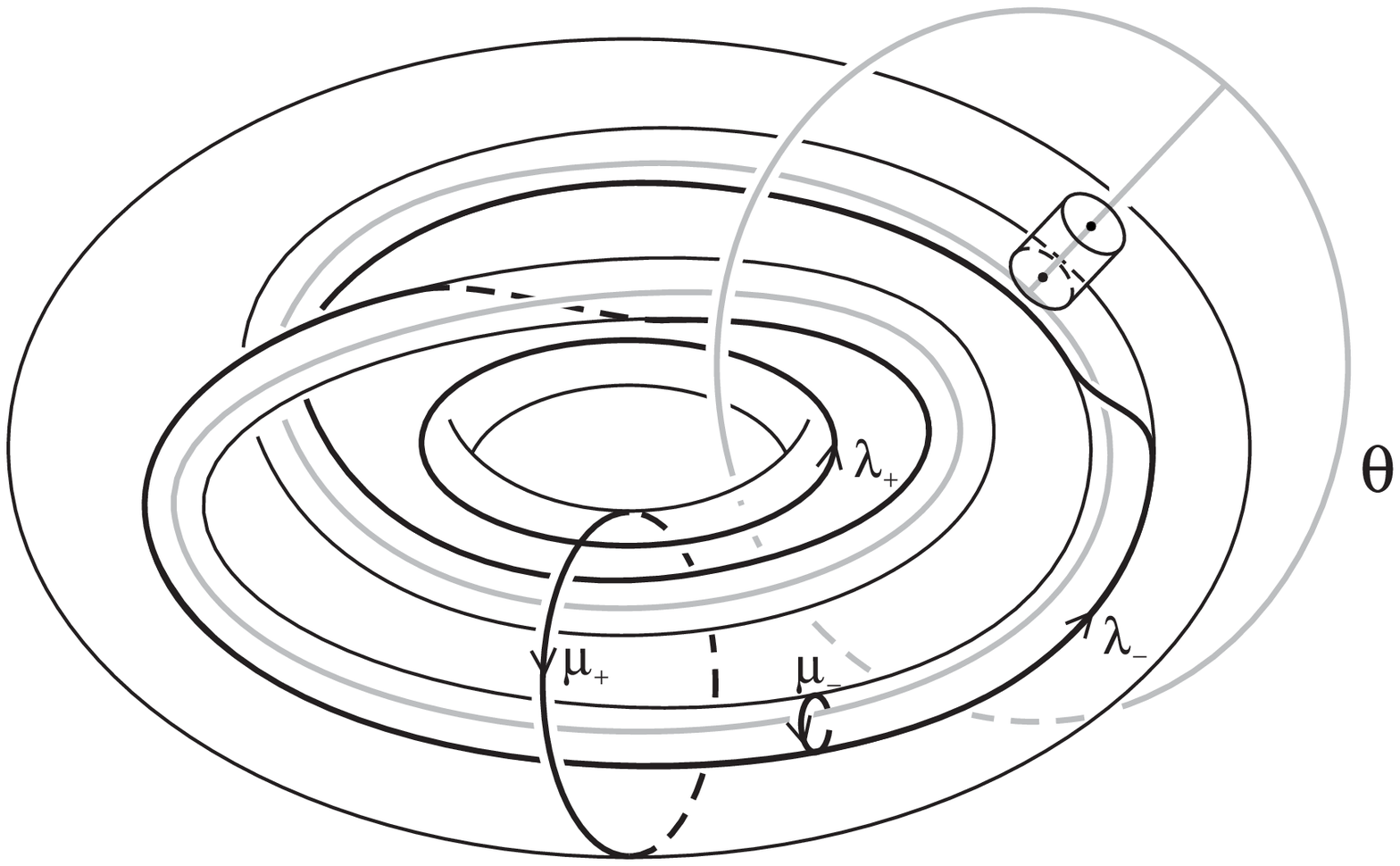}%
\begin{center}%
Figure 2
\end{center}
\end{center}

Since $f$ is a homeomorphism on $\mathbb S^3$ which sends the
unknotted solid torus $P$ to $f(P)$, we must have
$f(\lambda_+)=\lambda_-$, and $f(\mu_+)=\mu_-$. Now $M^*=X^*/f$ as
in Step~3 of \S 2.1.

Note that in Figure~2, $X^*$ is the complement of a graph $\Theta$
(shown  in gray in Figure 2) in $\mathbb S^3$, where $\Theta$
consists of the centerline of $f(P)$,  the centerline of $\mathbb
S^3\setminus P$, joined by the centerline $\gamma$  of the tunnel.

If we ignore the image of $X^*$ in Figure 2, but with $\Theta$,
$\lambda_\pm$ and $\mu_\pm$ remaining, then we have the Figure 3
below. Let $B^3$ be a 3-ball  containing the arc $\gamma$ in
$\Theta$, as indicated in Figure 3. It is an observation that the
complement of $\Theta$ is homeomorphic to the complement of two
unknotted arcs in the 3-ball $\mathbb S^3\setminus B^3$. Hence
$X^*$ is a handlebody of genus 2.

Two generators $a,b$ of $\pi_1(X^*)$ are indicated in Figure~3,
where we use the Wirtinger presentation \cite{R}, the base point
in $X^*$ being above the page. Representing $\lambda_{\pm}$ and
$\mu_{\pm}$ in terms of $a,b$, we have $\lambda_-=abab^{-1},\
\mu_-=b$; $\lambda_+=a,\ \mu_+=baba^{-1}$. By HNN extension, we
have
$$\pi_1(M^*)=<a,b,t|\ tat^{-1}=abab^{-1},\ tbaba^{-1}t^{-1}=b>,$$
and $\pi_1(\partial M^*)\cong{\mathbb Z}\oplus{\mathbb Z}$ is
generated by $t$ and $[\lambda_-,\mu_-]=[abab^{-1},b]$.
\begin{center}%
\includegraphics[totalheight=7cm]{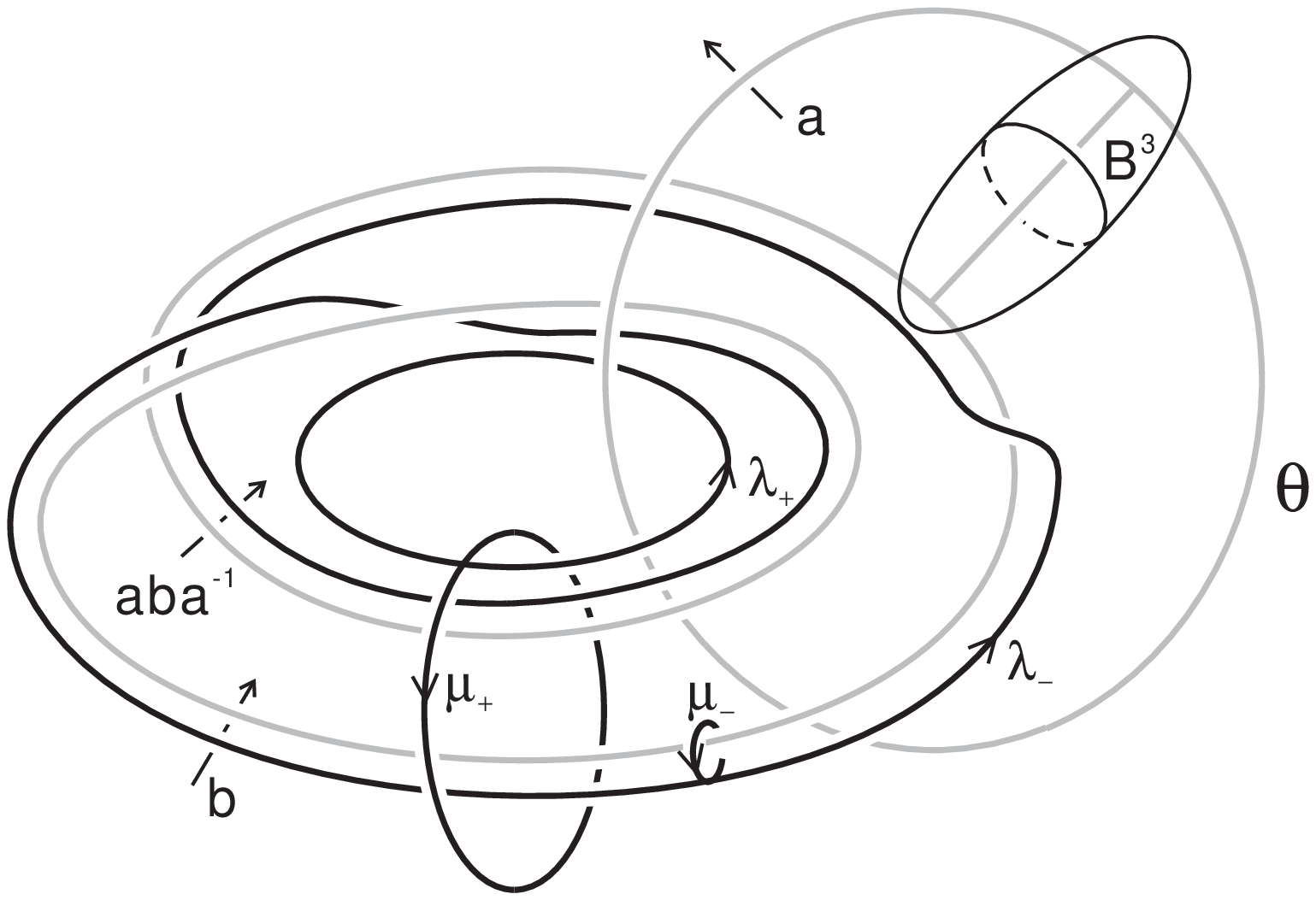}%
\begin{center}%
Figure 3
\end{center}
\end{center}
\vskip 0.3 true cm

Step 2. Compute $\pi_1(E)$ and $\pi_1(\partial E)\subset
\pi_1(E)$. We choose the projection of $9_{46}$ provided in [R, p.
211] rather than in the knot table of [R],  as Figure~4 below. The
Seifert surface $T$ of $9_{46}$ in Figure~4 is the 1-punctured
torus presented as a plumbing of two unknotted and untwisted bands
$B(\alpha)$ and $B(\beta)$ with oriented centerlines $\alpha$ and
$\beta$ respectively. $\pi_1(T)$ is generated by $\alpha$ and
$\beta$.
\begin{center}%
\includegraphics[totalheight=7cm]{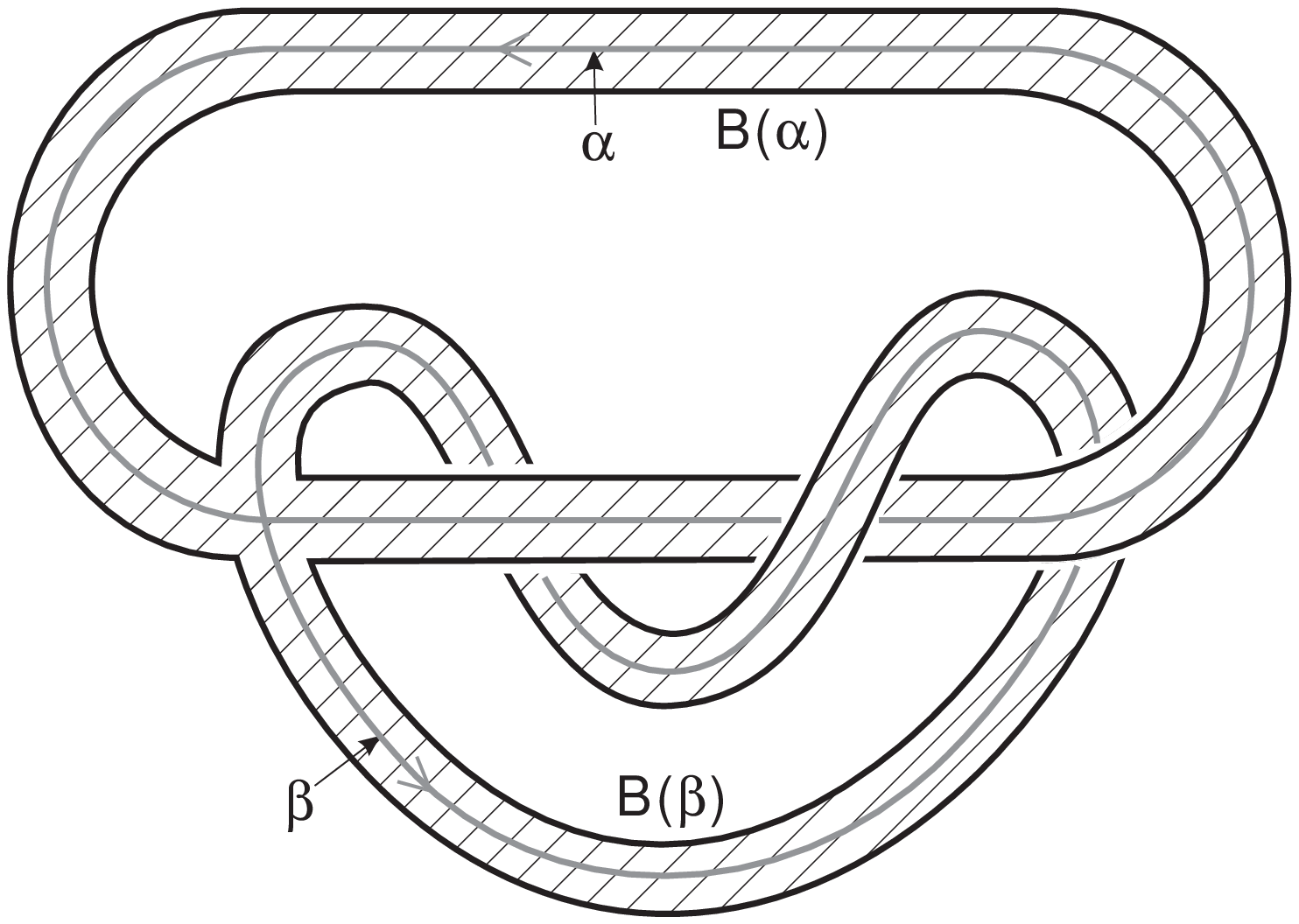}%
\begin{center}%
Figure 4
\end{center}
\end{center}
\vskip 0.3 true cm

Cutting $E$ open along $T$, we get a compact 3-manifold $Q$, which
is the complement of $T$ in $\mathbb S^3$, therefore $Q$ is also
homeomorphic to the complement of the one point union of the two
circles $\alpha\cup \beta$. By a handle sliding argument (see [R,
p. 95]) one can check that $Q$ is also a handlebody of genus 2.
Two generators $c,d$ of $\pi_1(Q)$ are indicated in Figure~5.

First pushing $\alpha$ and $\beta$ off $T$ towards the minus side
of $T$, we get two generators $\alpha_-,\beta_-$ of $\pi_1(T_-)$
in $\pi_1(Q)$; and then pushing $\alpha$ and $\beta$ off $T$
towards the plus side of $T$, we get two generators
$\alpha_+,\beta_+$ of $\pi_1(T_+)$ in $\pi_1(Q)$, all shown in
Figure 5. It can be easily computed that
$\alpha_-=cdcd^{-1},\beta_-=d,\alpha_+=c,\beta_+=dcdc^{-1}$. So
$$\pi_1(E)=<c,d,s|\ scs^{-1}=cdcd^{-1},\ sdcdc^{-1}s^{-1}=d>,$$
and $\pi_1(\partial E)\cong{\mathbb Z}\oplus{\mathbb Z}$ is
generated by $s$ and $[\alpha_-,\beta_-]=[cdcd^{-1},d]$.
\begin{center}%
\includegraphics[totalheight=7cm]{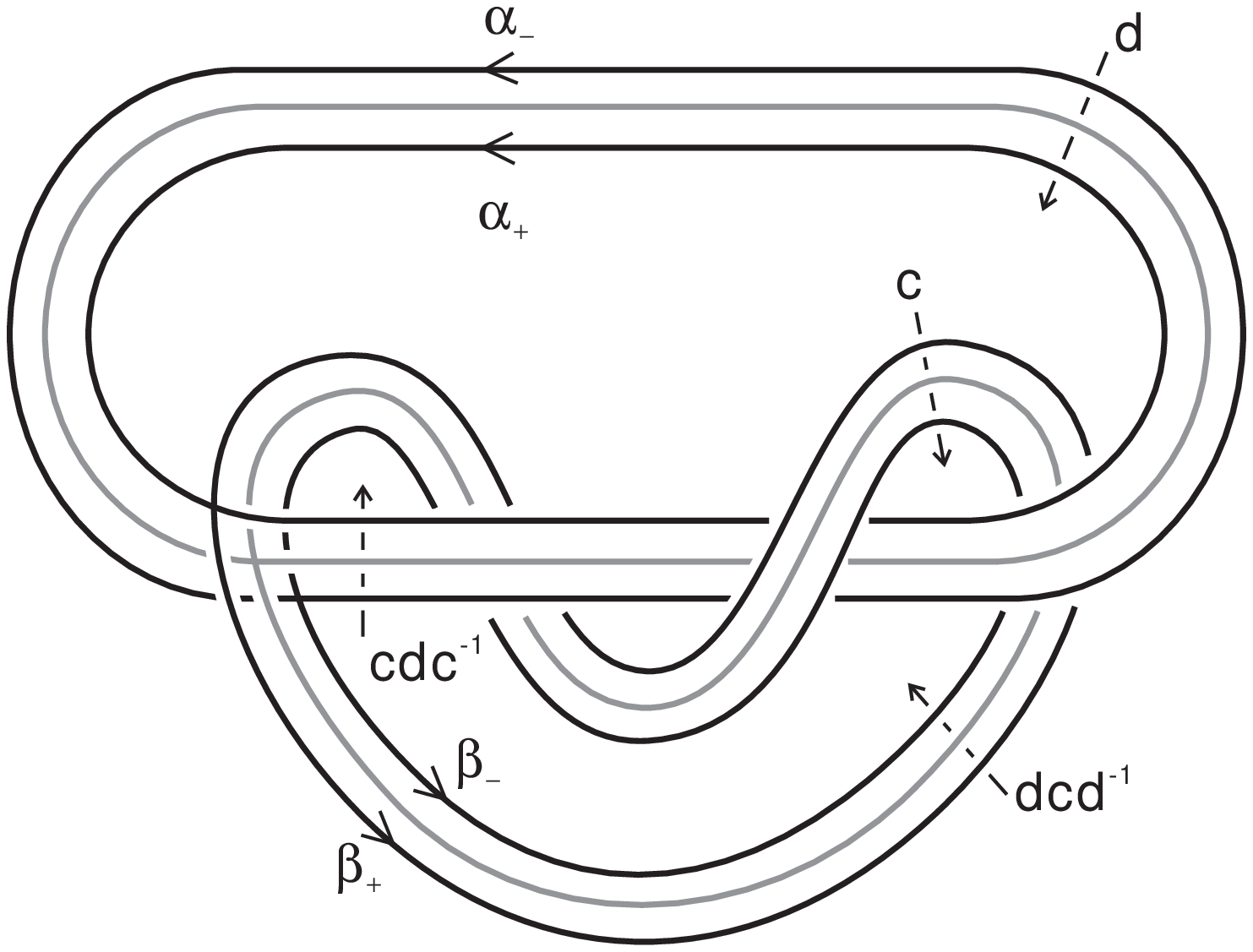}%
\begin{center}%
Figure 5
\end{center}
\end{center}
\vskip 0.3 true cm

Step 3. Now we have an isomorphism $$\phi:\pi_1(M^*)\to \pi_1(E),
\text{ such that } a\mapsto c, b\mapsto d,t\mapsto s,$$ which maps
$\pi_1(\partial M^*)$ isomorphically onto $\pi_1(\partial E)$.

Both $M^* ,E$ are $\mathbb P^2$-irreducible, sufficiently large
manifolds, so Waldhausen's theorem [H, Theorem 13.6] (or [BZ, p.
308 B7] more directly),  implies that $M^*$ is homeomorphic to
$E$. We finished the verification.


{\bf\S 2.3. Infinitely many genus 1 knots have Property $DIE$
}

Let $P$, $P'$, $X^*$ and $X^*/f$ be as given in Step 3 of \S 2.1.

{\bf Proposition 2.2}

{\it (1) If a meridian disk $D$ of $P$ meets the core of $P'$ in
exactly $2$ points tranversely, then $X^*/f$ is the complement of
a genus $1$ knot in a homotopy 3-sphere.

(2)  Furthermore if $X^*$ is homeomorphic to a handlebody of genus
$2$, then $X^*/f=E(k)$ for some genus $1$ knot $k\subset \mathbb
S^3$.

(3) There are infinitely many genus $1$ knots $k\subset\mathbb
S^3$ such that  $E(k)$ are obtained by the construction in \S
2.1.}

\begin{proof} Figure 6 indicates that  there are infinitely many
embeddings $P'\subset P$, such that both the conditions in
Proposition 2.2 (1) and (2) are satisfied. The verification of
$X^*$ to be the handlebody of genus 2 is the same as we did in
Figure 3 in \S 2.2. (Note that if we choose the tunnel jointing
$\partial P$ and $\partial P'$ to be knotted, then the condition
in (2) is not satisfied in general.)

\begin{center}%
\includegraphics[totalheight=7cm]{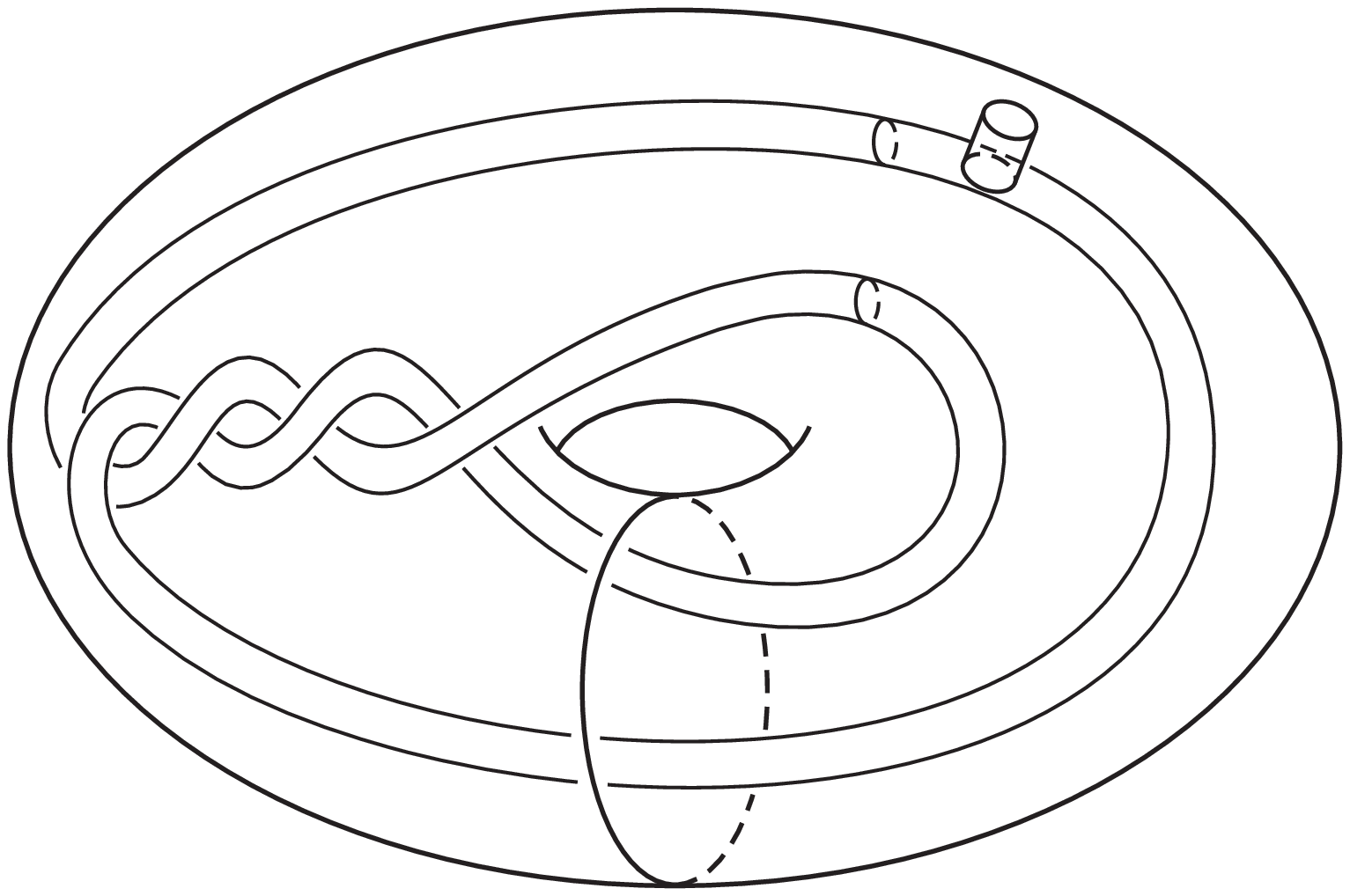}%
\begin{center}%
Figure 6
\end{center}
\end{center}

(1) We will find a presentation for $\pi_1(M^*)$ as in Step 1 of
\S 2.2 (but the process is simpler since we need less precise
information about the presentation.) First from Figure 6 we get
Figure 7, (as what we did from Figure 1 to Figure 3 in Step 1 of
\S 2.2,)  where $a,b,b'$ are elements in $G=\pi_1(X^*)$. Then as
in Step 1 of \S 2.2 we can compute $\pi_1(X^*/f)$ via HNN
extension as $$\pi_1(X^*/f)=<G,t|tat^{-1}=c,tbb't^{-1}=b>,$$ where
$c$ is the element in $G$ representing $\lambda_-$, and $t$ is
represented by a loop $\gamma$ in $\partial (X^*/f)=T^2$. Note
$\mu_+=bb'$ because the meridian disk intersects $P'$ twice.

\begin{center}%
\includegraphics[totalheight=7cm]{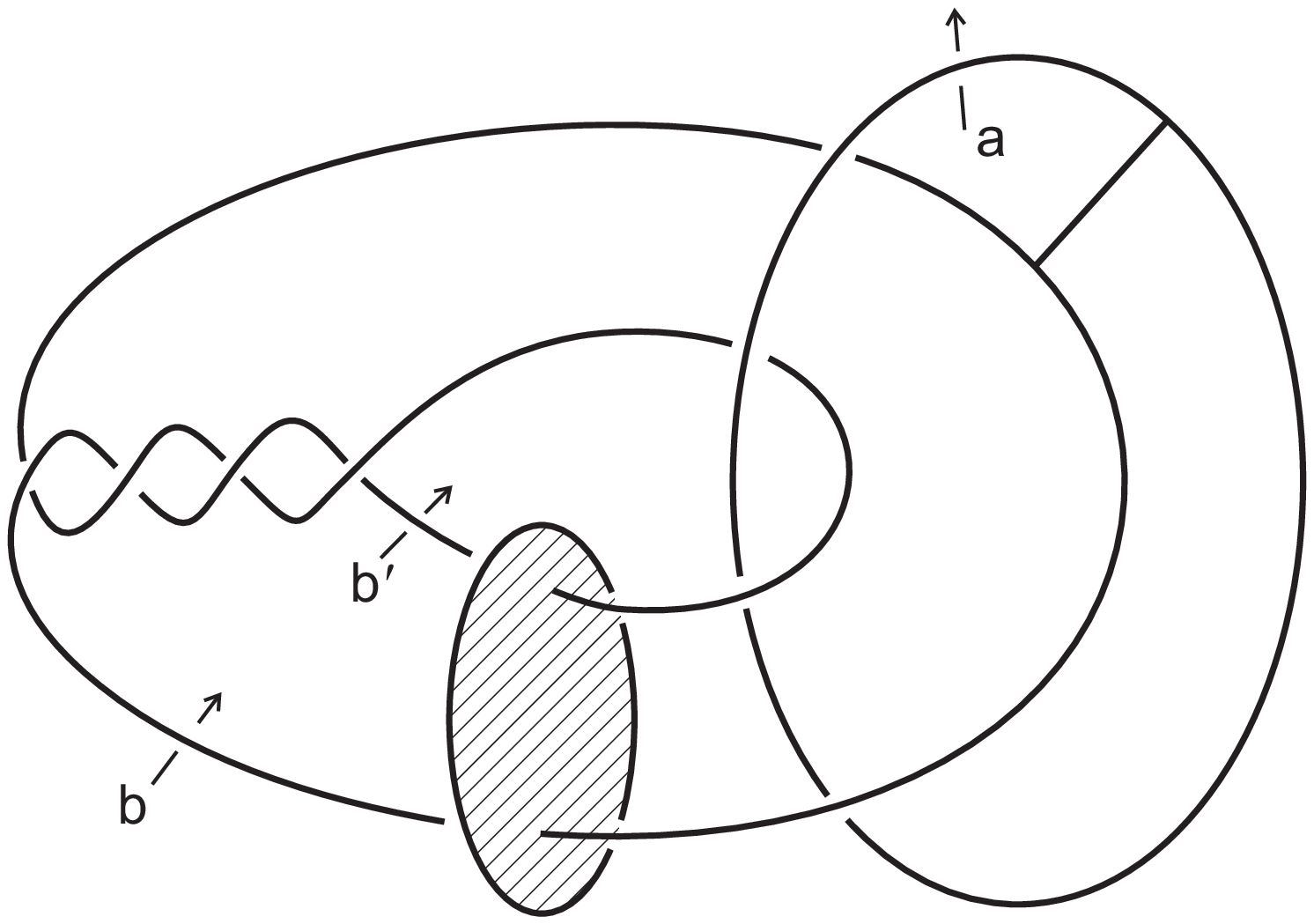}%
\begin{center}%
Figure 7
\end{center}
\end{center}

A Dehn filling along $\gamma$ will kill $t$ and provide a new
manifold $M_1=(X^*/f)(\gamma)$ with
$$\pi_1(M_1)=<G|a=c,bb'=b>\,=\,<G|a=c,b'=1>.$$

If we add a 2-handle to $X^*$ along the loop representing $b'$,
the new manifold is obviously a solid torus. So
$<G\,|b'=1>\cong\mathbb Z$. Thus $\pi_1(M_1)$ is a quotient group
of $\mathbb Z$. A computation in homology will show that
$\homo_1(M_1;\mathbb Z)=0$, hence $\pi_1(M_1)=1$. Thus $X^*/f$ is
the complement of a knot $k$ in the homotopy 3-sphere $M_1$.

(2) Furthermore suppose  $X^*$ is homeomorphic to a handlebody of
genus $2$. Note $X^*/f=X^*\cup_f N(\partial P\setminus
N(\Gamma))$, and $M_1=(X^*/f)(\gamma)$ can be viewed as a quotient
of $X^*\cup_f N(\partial P\setminus N(\Gamma))$ by identifying the
annulus $\partial (X^*/f)\cap N(\partial P\setminus N(\Gamma))$
with the annulus $\partial (X^*/f)\cap X^*$. Hence $M_1$ has a
Heegaard splitting $X^*\cup_h N(\partial P\setminus N(\Gamma))$ of
genus 2, where $h$ is determined by $f$ and $\gamma$. By Theorem 1
of \cite{BH} $M_1$ is a 2-fold cyclic covering of $S^3$, branched
over a 3-bridge link. It follows that $M_1$ is homeomorphic to
$\mathbb S^3$ by Thurston's orbifold theorem (see \cite{BP}), and
hence $X^*/f= E(k)$ for a knot $k$ in $\mathbb S^3$.

(3) We refine our notations related to Figure 6: Denote $P'$,
$X^*$ and $f$ by $P'_n$, $X^*_n$ and $f_n$, if the crossing number
of the core of $P'\subset P$ in Figure 6 is $n$, $n\in\mathbb Z$.
Then we have $X^*_n/f^*_n=E(k_n)$ for some knot $k_n\subset\mathbb
S^3$ according to (2). If there are only finitely many different
homeomorphism types for $E(k_n)$, then there are only finitely
many $E(k_n,0)$,  the zero surgery manifold on $k_n$. It follows
that the Gromov volumes $\{V(E(k_n, 0))\}$ take only finitely many
values. Note that $E(k_n,0)$ is homeomorphic to $(P\setminus
P'_n)/f_n$, and $E(k_n, 0)\setminus S_n = P\setminus P'_n$. Since
$\partial P$ is incompressible in $P\setminus P'_n$,
$V(E(k_n,0))=V(P\setminus P'_n)$ by a theorem of Soma [S, Theorem
1], it follows that $\{V(P\setminus P'_n)\}$ take only finitely
many values.

\begin{center}%
\includegraphics[totalheight=4.3cm]{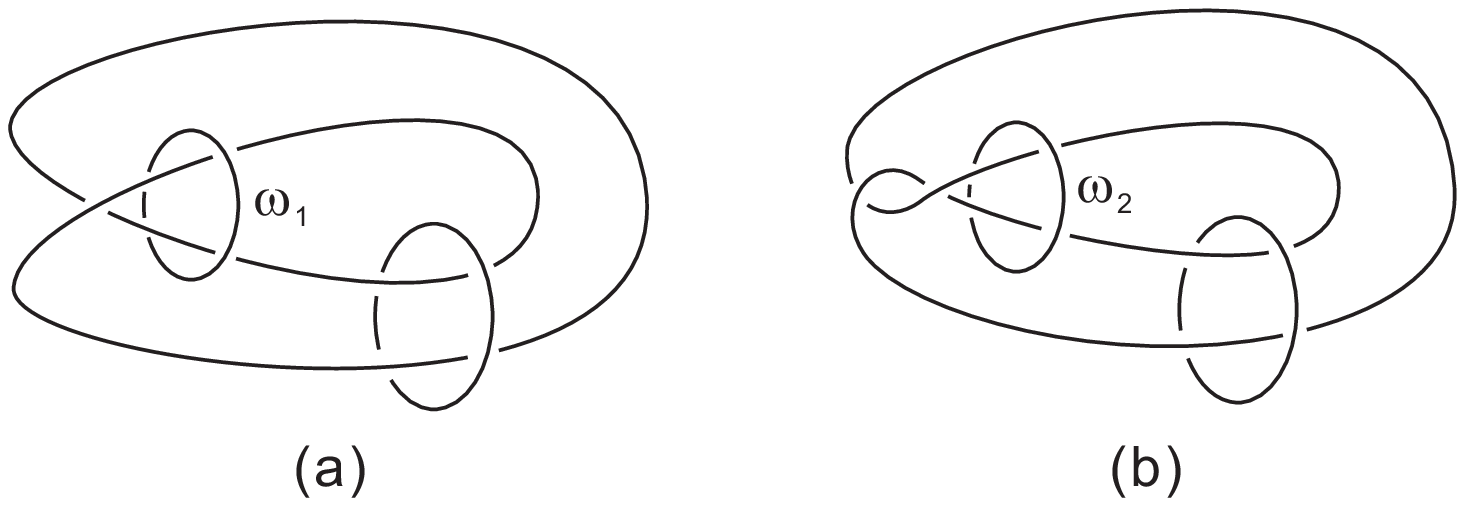}%
\begin{center}%
Figure 8
\end{center}
\end{center}

Consider two 3-component links $L_1$ and $L_2$ with marked
components $\omega_1$ and $\omega_2$ respectively, indicated in
Figure 8, (a) and (b). Note $\omega_i$ is unknotted in $\mathbb
S^3$, the standard arguments (see [R, Chap. 9]) show that
$$P\setminus P'_{2n+1}=E(L_1)(\omega_1, 1/n),
\qquad P\setminus P'_{2n+2}=E(L_2)(\omega_2, 1/n),$$ where
$E(L_i)(\omega_i, 1/n)$ is the $1/n$-Dehn filling along
$\omega_i$. It is also known that both $L_1$ and $L_2$ are
hyperbolic links. (This fact can be checked by SnapPea \cite{W}.)
According to Thurston's theory about Gromov volume on 3-manifolds
(see [Th, Chapters 5 and 6]), we have

(i) $V(E(L_i)(w_i,1/n))< V(E(L_i))$,

(ii) $\lim_{n\to \infty}V(E(L_i)(w_i,1/n))= V(E(L_i))$.

It follows that $\{V(P\setminus P'_n)\}$ take infinitely many
values, a contradiction.\end{proof}

{\bf Remark 1.} All knots constructed in Proposition 2.2 bound the
gunus 1 surface. All knots $k_n$ in Propositon 2.2 (3) are
non-fibred, see the end of \S 3, and also $k_1=9_{46}$.

\vskip 0.2 true cm

\centerline {\bf \S 3. ``Most" Genus 1 knots do not have Property
$IE$}

Let $k$ be a non-fibred knot of genus 1 in $\mathbb S^3$. Recall
the notations $E(k)$, $S$, $X=E(k)\setminus S, \widetilde E(k)$
defined in the beginning of the paper. Suppose also $S$ is of
genus 1.

Let $S_n$ ($n=\dots,-2,-1,0,1,2,\dots$) denote the copies of $ S$
in $\widetilde E(k)$. For integers $m<n$, let $X_{[m,n]}$ denote
the sub-manifold of $\widetilde E(k)$ between $ S_m$ and $ S_n$,
and $ A_{[m,n]}$ denote the annulus bounded by $\partial
S_m\sqcup\partial S_n$ on $\partial X_{[m,n]}$. Assume $\widetilde
E(k)$ is already embedded in $\mathbb S^3$, and $Y_{]m,n[}=\mathbb
S^3\setminus X_{[m,n]}$. We always use $X_n$ to denote
$X_{[n,n+1]}$ for simplicity. The readers should be aware that the
subscript $n$ here has different meaning from the $n$ in the last
section.

{\bf Lemma 3.1} {\it For any integer $N>0$, $\partial S_0$ bounds
a disk $ D$ in $Y_{]-N,N[}$.}
\begin{proof}
Consider the separating surfaces $ S^*_{n}= A_{[N,n]}\cup S_{n}$
($n=N+1,N+2,\dots$) in $Y_{]-N,N[}$. They are mutually
non-parallel, since $k$ is non-fibred. Since each $S_n^*$ has the
first Betti number 2, Freedman-Freedman's version of Kneser-Haken
finiteness theorem \cite{FF} implies that they must be
compressible in $Y_{]-N,N[}$ when $n$ is sufficiently large.
Suppose $D$ is a compressing disk of $ S^*_{n}$. If $\partial D$
is parallel to $\partial S^*_{n}$ on $ S^*_{n}$, then the lemma is
proved, since $\partial S^*_{n}$ is parallel to $\partial S_0$ on
$\partial Y_{]-N,N[}$. If $\partial D$ is not parallel to
$\partial S^*_{n}$, surger $ S^*_{n}$ along $\partial D$, we still
get a disk $D'$ in $Y_{]-N,N[}$, with $\partial D'=\partial
S^*_n$, since $S_n^*$ is a 1-punctured torus.
\end{proof}

Now fix an $N$ sufficiently large, we can thicken $ D\cup
A_{[-N,N]}$ in $Y_{]-N,N[}$ to get a 2-handle $ D\times I$, which
is attached to $X_{[-N,N]}$ along the annulus $ A_{[-N,N]}$. Let $
D_{-N},\dots, D_{N}$ be a collection of $ D\times\{t\}$'s in the
2-handle, so that $\partial D_i=\partial S_i,\; i=-N,\dots,N$.
From now on, all subscripts in this section are bounded by $N$, as
is understood.

Let $\widehat{S}_i$ denote the torus $ S_i\cup D_i$. Let $\widehat
X_i$ be the manifold bounded by $\widehat{S}_i$ and $\widehat{
S}_{i+1}$ in $\mathbb S^3$, and more generally, $\widehat
X_{[m,n]}$ be the manifold bounded by $\widehat{ S}_m$ and
$\widehat{ S}_n$ in $\mathbb S^3$.

{\bf Lemma 3.2} {\it $\widehat X_i$ is irreducible and
$\partial$-irreducible. Moreover $\widehat X_i$ is not a product.}

\begin{proof} Since $S$ is a minimal genus Seifert surface of $k$,
 $E(k)$ admits a taut foliation $\mathcal F$
 such that $S$ is a leaf of $\mathcal F$ and $\mathcal F |
\partial E(k)$ is foliated by circles by [G1, Theorem 3.1]. Then $\mathcal F$ can
be extended to a taut foliation $\widehat{\mathcal F}$ on
$E(k,0)$, the zero surgery manifold on $k$, such that $\widehat S$
is a leaf of $\widehat{\mathcal  F}$, where $\widehat S$ is
obtained by capping disc on $S$. Moreover since $E(k)$ is not
fibred, $E(k,0)$ is not fibred by [G1, Corollary 8.19], in
particular $E(k,0)\ne S^2\times S^1$. By Novikov's theorem [N],
each leaf of the taut foliation $\widehat{\mathcal F}$ is
$\pi_1$-injective in $E(k,0)$ and $\pi_2($E(k,0)$)=0$. Then
$E(k,0)$ is irreducible by the sphere theroem [H, Chap. 3], and
furthermore $\widehat S$ is incompressible. It follows that
$E(k,0)\setminus \widehat S$ is irreducible,
$\partial$-irreducible, and is not a product.

Since each $\widehat X_i$ is homeomorphic to $E(k,0)\setminus
\widehat S$, Lemma 3.2 is proved.
\end{proof}

Each $\widehat{ S}_i$ separates $\mathbb S^3$ into 2 components.
We say the component containing $\widehat X_i$ {\sl lies on the
plus side of $\widehat{ S}_i$}, the component containing $\widehat
X_{i-1}$ {\sl lies on the minus side of $\widehat{ S}_i$}.
$\widehat{ S}_i$ bounds a solid torus on the plus side or the
minus side, since every torus in $\mathbb S^3$ bounds a solid
torus. In fact, we can prove the stronger

{\bf Proposition 3.3} {\it Each $\widehat{ S}_i$ bounds solid tori
on both sides.}
\begin{proof}
Without loss of generality, we can assume $\widehat{ S}_0$ bounds
a solid torus $P_0$ on the minus side. Our argument proceeds in
the following steps.

{\bf Step 1.}\; For each $n<0$, $\widehat{ S}_n$ bounds a solid
torus $P_n$ on the minus side.

Otherwise, assume some $\widehat{ S}_n$ does not bounds a solid
torus on the minus side, then $\widehat{ S}_n$ bounds a solid
torus $P_+$ on the positive side. Hence $\widehat{ S}_n$ cuts
$P_0$ into 2 parts: $\widehat X_{[n,0]}$ and $P_0\setminus\widehat
X_{[n,0]}=\mathbb S^3\setminus P_+$. By Lemma 3.2, $\widehat S_n$
is incompressible in $\widehat X_{[n,0]}$; $\widehat S_n$ is also
incompressible in $\mathbb S^3\setminus P_+$ since $P_+$ is
knotted. So $P_0=(\mathbb S^3\setminus P_+)\cup_{\widehat{
S}_n}\widehat X_{[n,0]}$ cannot have $\pi_1=\mathbb Z$.

By Step 1, we have a nested sequence of solid tori
$$...\subset P_{n-1}\subset P_n \subset P_{n+1} \subset ...\subset P_0.$$
We assume that these tori adapt the orientation of $\mathbb S^3$.
Let $\mu_n,\, \lambda_n \subset \widehat{ S}_n$ be a oriented
meridian-longitude system of $P_n$, $n< 0$,  so that

(1) the algebraic intersection number of $\mu_n$ and $\lambda_n$
is 1,

(2) the linking number of $\lambda_{n}$ and $\mu_{n+1}$, which is
defined as the winding number of $P_n$ in $P_{n+1}$,  is $\ge 0$.

{\bf Step 2.}\; Suppose $P_n$ has winding number $w_n$ in
$P_{n+1}$, $n<0$. Then all $w_n$ are equal, denoted by $w$.

Clearly $P_n\setminus P_{n-1}$ is homeomorphic to the complement
in $ \mathbb S^3$ of a 2-component link with linking number
$w_{n-1}$, so $H_1(P_n\setminus P_{n-1}; \mathbb Z)$ has a basis
$\lambda_n, \mu_{n-1}$, and $H_1(P_n\setminus P_{n-1}, \widehat
S_n; \mathbb Z)$ is isomorphic to $Z_{w_{n-1}}$.

Note that the deck translation $\tau:\widetilde E(k)\to\widetilde
E(k)$, which sends $X_i$ to $X_{i+1}$, induces a homeomorphism
$\widehat{\tau}: \widehat X_{n-1}=P_n\setminus P_{n-1}\to \widehat
X_n=P_{n+1}\setminus P_n$ with $\widehat \tau_n(\widehat S_n)=
\widehat S_{n+1}$ for each $n<0$. It follows that $w_n=w_{n-1}$.

{\bf Step 3.}\;  We claim that $\widehat{\tau}$ sends $\mu_{n-1}$
to $\mu_{n}$ for $n\le -1$. There are 2 cases:

Case 1. $w=0,1$. Now $P_n$ can not be a braid in $P_{n+1}$,
otherwise $w=1$ and $\widehat X_n$ is a product $T^2\times I$,
contradicts to Lemma 3.2. Then the  results in \cite{G2} imply
that only trivial surgery on $P_n$ yields solid torus. Since the
Dehn surgery on the knot $P_{n}$ in the solid torus $P_{n+1}$
along $\widehat \tau(\mu_{n-1})$ again yields solid torus,
$\widehat{\tau}(\mu_{n-1})=\mu_{n}$ for $n\le-1$.

Case 2. $w\ge 2$. Fix $n<0$. Now $\lambda_i,\mu_i$ is a basis of
$\homo_1(\widehat{ S}_i;\mathbb \mathbb Z)$, $i\le0$.
$$\widehat{\tau}_*(\lambda_n)=p\lambda_{n+1}+q\mu_{n+1},\qquad
\widehat{\tau}_*(\mu_n)=r\lambda_{n+1}+s\mu_{n+1},\qquad
ps-qr=1.$$

For each integer $m>0$, since $\mu_n$ is a $w^m$ multiple in
$\homo_1(\widehat X_{[n-m,n]};\mathbb Z)$,
$\widehat{\tau}_*(\mu_n)$ is also a $w^m$ multiple in
$\homo_1(\widehat X_{[n-m+1,n+1]};\mathbb Z)$. Since $\mu_{n+1}$
is already a $w^m$ multiple in $\homo_1(\widehat
X_{[n-m+1,n+1]};\mathbb Z)$, $r\lambda_{n+1}$ is also a $w^m$
multiple.

Since $\{\lambda_{n+1}, \mu_{n-m+1}\}$ is a basis of
$\homo_1(\widehat X_{[n-m+1,n+1]})$ for $m>0$, $r$ should be a
$w^m$ multiple. Since $r$ is a given integer, letting $m$ be
sufficiently large, we must have $r=0$. Then $p=s=\pm1$, i.e.,
$\widehat{\tau}_*(\mu_n)=\pm \mu_{n+1}$, the conclusion holds.

{\bf Step 4.}\; When $n>0$, $\widehat{ S}_n$ bounds a solid torus
on the minus side.

There is a properly embedded planar surface $G$ in $\widehat
X_{-2}$, $G\cap\widehat{ S}_{-1}=\mu_{-1}$, $G\cap\widehat{
S}_{-2}$ consists of parallel copies of $\mu_{-2}$. By Step 3,
$\widehat{\tau}(G)$ is a planar surface in $\widehat X_{-1}$,
$\widehat{\tau}(G)\cap\widehat{ S}_{-1}$ consists of parallel
copies of $\mu_{-1}$. $\widehat{\tau}(G)\cap\widehat{ S}_{0}$
bounds a disk on the minus side of $\widehat{ S}_0$, since each
copy of $\mu_{-1}$ bounds a disk in $P_{-1}$. So $\widehat \tau
(\mu_{-1}) =\widehat{\tau}(G)\cap\widehat{ S}_{0}=\mu_0$. Let
$\mu_n=\widehat{\tau}^n(\mu_0)$ for $n>0$, the same argument as
above shows that $\mu_n$ bounds a disk on the minus side of
$\widehat{ S}_n$, by induction.

{\bf Step 5.}\; All $\widehat{ S}_n$ bounds solid tori on both
sides, $n\in\mathbb N$.

By Lemma 3.2, $\widehat S_n$ and $\widehat S_m$ are not parallel
for $m\ne n$. By Haken's finiteness theorem, $\widehat{ S}_n$ is
compressible in $\mathbb S^3\setminus P_0$ when $n$ is
sufficiently large. The compressing disk can not lie on the minus
side, since $\widehat X_{[0,n]}$ is $\partial$-irreducible by
Lemma 3.2. So $\widehat S_n$ bounds a solid torus on the plus side
when $n$ is sufficiently large. Now proceed from Step 1 to Step 4,
but reverse the direction, to get our conclusion.
\end{proof}

{\bf Theorem 3.4} {\it Suppose $k$ is a non-fibred knot of genus
$1$ in $\mathbb S^3$. If $k$ have Property $IE$, then $k$ has
Property $DIE$. Indeed, $E(k)$ can be  obtained by the
construction in \S 2.1.

Moreover, the winding number $w$ involved is either $0$ or $2$.
Correspondingly, the Alexander invariant of $k$ is either 0 or
$\mathbb Z[t,t^{-1}]/(2t-1)\oplus\mathbb Z[t,t^{-1}]/(t-2)$, and
the Alexander polynomial of $k$ is either $1$ or $2t^2-5t+2$.}

\begin{proof}
Suppose $\widetilde E(k)$ is embedded into $\mathbb S^3$. We keep
the notations in the proof of Proposition 3.3. First, extend
$\tau|X_0: (X_0, S_0) \to (X_1, S_1)$ to a homeomorphism $\widehat
\tau_1: (\widehat X_0, \widehat S_0)\to (\widehat X_1,\widehat
S_1)$ as in the proof of Proposition 3.3.

 According to Proposition 3.3, each
$\widehat{ S}_n$ bounds a solid torus $P^-_n$ on the minus side,
and a solid torus $P^+_n$ on the plus side. Suppose
$\mu_n^-,\mu_n^+\subset S_n\subset \widehat S_n$ are meridians of
$P^-_n,P^+_n$ respectively. By Step 3 (and its counterpart in Step
5) of Proposition 3.3,
$$\widehat{\tau}(\mu_n^-)=\mu_{n+1}^-,\qquad \widehat{\tau}
(\mu_n^+)=\mu_{n+1}^+.\eqno(3.1)$$

Hence we can further extend $\widehat \tau_1$ to $\widehat \tau_2
: P_0^+\to P_1^+$, and finally we extend $\widehat \tau_2$ to
$f:\mathbb S^3\to\mathbb S^3$ since both $P_0^+$ and $P_1^+$ are
unknotted. Now we can reconstruct $E(k)$ from $f$ as in \S 2.1, so
$k$ has Property $DIE$. We have finished the proof of the first
part of Theorem 3.4.

By Step 2 (and its counterpart in Step 5) of Proposition 3.3, the
winding number of $P^-_{n}$ in $P^-_{n+1}$ is a constant $w^-$,
and the winding number of $P^+_{n+1}$ in $P^+_{n}$ is a constant
$w^+$. It is easy to see that  both $w^-$ and $w^+$ are the
linking number between $\mu_{n+1}^-$ and $\mu_{n}^+$ (see the
paragraph after Step 1 in the proof of Proposition 3.3), we have
$w^-=w^+=w$. Since $\widehat \tau| : \widehat S_n\to \widehat
S_{n+1}$ is orientation preserving, by (3.1) we have
$$\widehat{\tau}^{-1}_*([\mu_n^+])=\pm w[\mu_{n}^+],\qquad
\widehat{\tau}_*([\mu_n^-])=\pm w[\mu_{n}^-],\eqno(3.2).$$

Note that $X_n\hookrightarrow\widehat X_n$ induces an isomorphism
on 1-dimensional homology. Then by (3.2) the  Alexander invariant
of $k$ has presentation [R, Chap 7]
$$\homo_1(\widetilde E(k);\mathbb Z[t,t^{-1}]) =<\mu^+_n, \mu^-_n, t |\, t^{-1}([\mu_n^+])=\pm
w[\mu_{n}^+],\, t([\mu_n^-])=\pm w[\mu_{n}^-]>,$$ and the
Alexander matrix of $k$ is
$$\left(\begin{array}{cc}wt\mp1&0\\
0&t\mp w
\end{array}\right).
$$

Since $\Delta_k(1)=\pm1$, $w$  can only be $0$ or $2$, and the
corresponding Alexander polynomials are $1$ or $2t^2-5t+2$
respectively, and the Alexander invariant of $k$ are either $0$ or
$\mathbb Z[t,t^{-1}]/(2t-1)\oplus\mathbb Z[t,t^{-1}]/(t-2)$. We
have finished the proof of Theorem 3.4.
\end{proof}

{\bf Corollary 3.5} {\it Among all genus $1$ non-fibred knots in
$\mathbb S^3$,

(1) up to ten crossings, $9_{46}$ is the only one that has
Property $IE$,

(2) no alternating knot has Property $IE$.}

\begin{proof} (1) For knots with $\le10$ crossings, no non-fibred
knot has Alexander polynomial $1$, and only  $6_1$ and $9_{46}$
have Alexander polynomial $2t^2-5t+2$, see the tables in [BZ] and
in [R]. But their Alexander invariants are not isomorphic (see [R,
p.211]), so $6_1$ does not have Property $IE$.  Then by \S 2.2 (1)
follows.

(2) If a genus 1 non-fibred knot $k$ has Property $IE$, then
$\Delta_k(-1)=1\quad \text{or}\quad 9$. Now suppose $k$ is
alternating, by a theorem of R.H. Crowell, (see [BZ, Proposition
13.30]) $\Delta_k(-1)$ is not smaller than the crossing number of
$k$, and $9_{46}$ is not alternating. Hence (2) follows from (1).
\end{proof}

Recall the two infinite families of knots $k_{2n}$ and $k_{2n+1}$
with Property $IE$, as well as the notion $P'_n$, defined in the
proof of Proposition 2.2 (3). Since the winding number of
$P'_{2n}$ is 0 and the winding number of $P'_{2n+1}$ is 2,
according to the calculation in the proof of Theorem 3.4 we have
$\Delta_{k_{2n}}(t)=1$ and $\Delta_{k_{2n+1}}(t)=2t^2-5t+2$.

{\bf Corollary 3.6} {\it Among non-fibered genu 1 knots, both the
subsets defined by $\Delta_{k}(t)=1$ and by
$\Delta_{k}(t)=2t^2-5t+2$ have infinitely many elements with
Property $IE$.\hfill\qedsymbol}

\begin{center}{\bf \S 4. A remark on connected sums}\end{center}

{\bf Lemma 4.1} {\it Suppose $k_1$ and $k_2$ are two knots in
$\mathbb S^3$.

(1) If $k_1\#k_2$ has Property $IE$, then both $k_1$ and $k_2$
have Property $IE$.

(2) If $k_1$ has Property $IE$ and $k_2$ is fibred, then
$k_1\#k_2$ has Property $IE$.}

Note that there are fibred knots of any genus (just consider the
connected sum of genus 1 fibred knots), and that $k_1\# k_2$ is
fibred if and only if both $k_1$ and $k_2$ are fibred (follows
from the definitions of connected sum, fibred knot, and Stallings'
fibration Theorem [H, Theorem 11.1]). Then by the main results in
\S 2, \S3 and Lemma 4.1 we have the following

{\bf Corollary 4.2} {\it Among non-fibred knots of genus $g$ for
any given integer $g>0$, both the subsets defined by  having
Property $IE$ and not having Property $IE$ have infinitely many
elements.\hfill\qedsymbol}

\vskip 0.2 true cm

\noindent{\it Proof of Lemma 4.1.} Denote $E(k_i)$ by $E_i$. Let
$N_i=N(\mu_i)$ be the regular neighborhood of the meridian
$\mu_i\subset \partial E_i$ in $E_i$. Let $E^*_i=E_i\setminus
N_i$, and $A_i=E_i^*\cap N_i$. Then $E^*_i$ is homeomorphic to
$E_i$ and $A_i$ is an annulus. By definition of the connected sum,
we have $E(k_1\# k_2)=E^*_1\cup_h E^*_2$, where $h$ is a
homeomorphism identifying $A_1$ and $A_2$.

\begin{center}%
\includegraphics[totalheight=8cm]{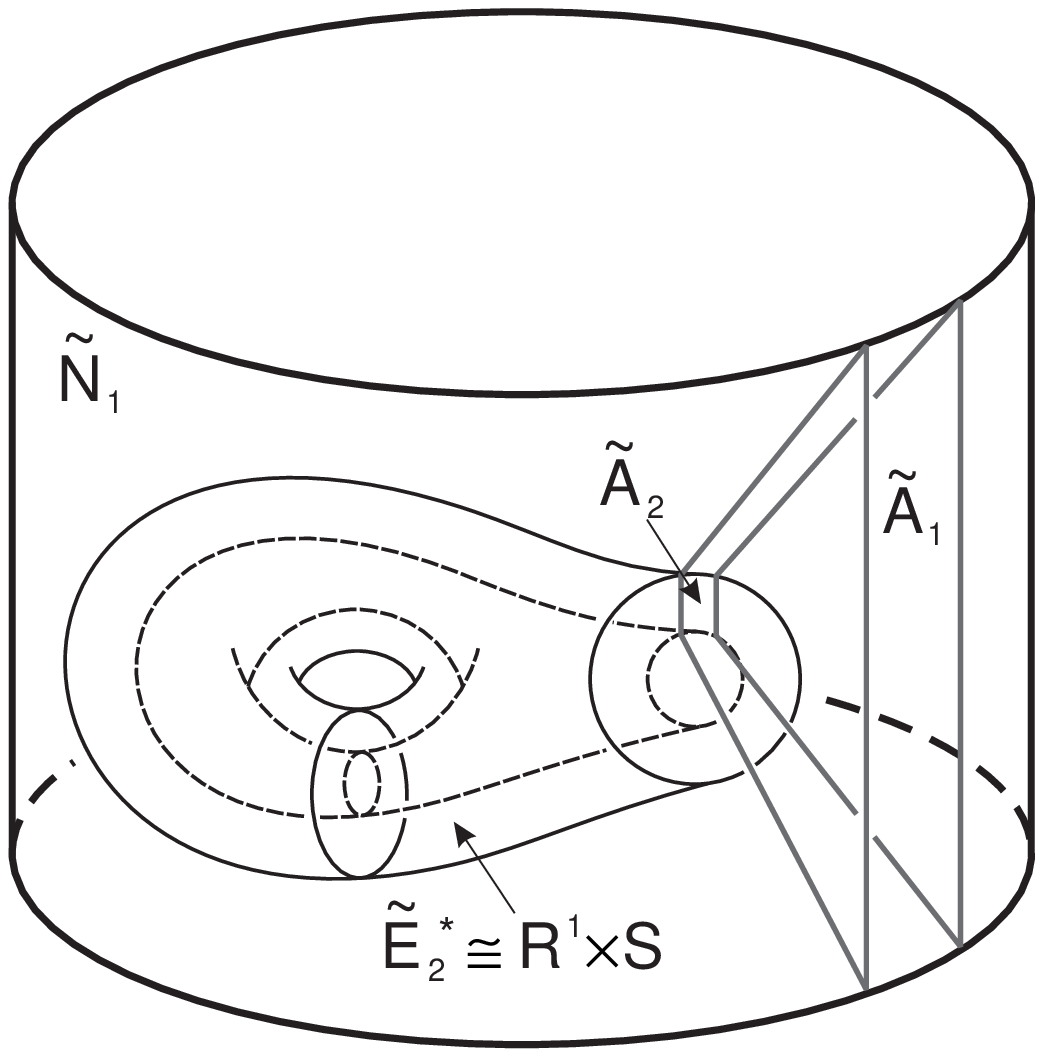}%
\begin{center}%
Figure 9
\end{center}
\end{center}

Let $p_i: \widetilde E_i\to E_i$ be the infinite cyclic covering,
and let $\widetilde E^*_i, \widetilde N_i, \widetilde A_i$ be the
preimage of $E^*_i, N_i, A_i$ under $p_i$. Clearly the restriction
of
$$p_i : (\widetilde E_i, \widetilde E^*_i, \widetilde N_i, \widetilde A_i)\to (E_i, E^*_i, N_i,
A_i)$$ is the infinite cyclic covering on each of the four
corresponding pairs. Moreover $\widetilde E^*_i, \widetilde N_i,
\widetilde A_i$ are homeomorphic to $\widetilde E_i, R^1\times
D^2, R^1\times I$ respectively and $\widetilde E(k_1\#
k_2)=\widetilde E_1^*\cup_{\widetilde h} \widetilde E_2^*$, where
$\widetilde h$ is a homeomorphism identifying $\widetilde A_1$
with $\widetilde A_2$. Hence (1) follows.

We are going to prove (2). Now $\widetilde E^*_2= R^1\times S$ for
a once punctured surface $S$ with $\widetilde A_2=R^1\times I$
properly embedded in $R^1\times \partial S$.

Since there is an embedding $e: \widetilde E^*_2=R^1\times S \to
\widetilde N_1$ such that $e$ sends $\widetilde A_2$ to
$\widetilde A_1\subset
\partial \widetilde N_1$ homeomorphically, and $e(\widetilde E^*_2)\cap \widetilde N_1=\widetilde
A_1$ (see Figure 9), $\widetilde E(k_1\# k_2)=\widetilde
E_1\cup_{\widetilde h} \widetilde E_2$ can be embedded into
$E^*_1\cup\widetilde N_1=\widetilde E_1$. Hence (2) follows.
\hfill\qedsymbol

\begin{center}{\bf \S 5. A partial negative answer to Question 2}\end{center}

 In this section we use the notations in the first two paragraphs
of \S 1. We will use $\homo_i(\cdot)$ to denote
$\homo_i(\cdot\;;\mathbb Q)$. Recall the following standard fact:
Let $$\dots\to A \to B \to  C\to \cdots$$ be an exact sequence of
vector spaces. Then
$$\dim A + \dim C\ge \dim B.
 \eqno (*)$$

{\bf Theorem 5.1} {\it Suppose $M$ is a compact 3-manifold, $S$ is
a connected non-separating 2-sided proper surface in $M$. Let
$X=M\setminus S$.

(1) In the case $\partial M\ne \emptyset$, if $[S\cap T]\ne 0\in
\homo_1(\partial M;\mathbb Z)$ for each boundary component $T$ of
$M$ and $\beta_1(X)>\beta_1(S)-\chi(\partial M),$ then $\widetilde
M_S$ cannot be embedded into any compact 3-manifold.

(2) In the case $\partial M= \emptyset$, if
$\beta_1(X)>\beta_1(S),$ then $\widetilde M_S$ cannot be embedded
into any compact 3-manifold.}

\begin{proof} Suppose $\partial M\ne\emptyset$, $\widetilde M=\cup_{k=-\infty}^{+\infty}X_k$ can be
embedded into a compact $3$-manifold $Y$. We may assume $\partial
Y=\emptyset$. Denote $\cup_{k=1}^mX_k$ by $P_m$.

We need first estimate $\beta_1(P_m)$. From $P_m=P_{m-1}\cup
X_{m}$ and $S_m=P_{m-1}\cap X_{m}$, we have the Mayer-Vietoris
sequence:
$$
\cdots{\to}\homo_1(S_m)
\to\homo_1(P_{m-1})\oplus\homo_1(X_{m})\to\homo_1(P_m) \to\cdots.
$$

By ($*$), we have the inequality:
$$\beta_1(P_m)\ge\beta_1(P_{m-1})+\beta_1(X)-\beta_1(S).$$

Hence we easily deduce:
$$\beta_1(P_m)\ge m\beta_1(X)-(m-1)\beta_1(S). \eqno (1)$$

We need then estimate $\beta_1(\partial P_m)$.

Cutting $\partial M$ open along $\partial S$, we get a surface
$T'$. $\partial P_m$ is the union of $S^-_1\sqcup S^+_m$ and $m$
copies of $T'$.  Note that the cutting and gluing of surfaces are
all along circles, which have Euler characteristic 0. So
\begin{eqnarray*}
\chi(\partial P_m) &= &\chi(S^-_1\sqcup S^+_m)+m\chi(T')\\
&= &2\chi(S)+m\chi(\partial M)\\
&= &2(1-\beta_1(S))+m\chi(\partial M)\\
\end{eqnarray*}

Then one can verify that
$$\beta_1(\partial P_m)=2\beta_0(\partial P_m)-\chi(\partial P_m)
=2\beta_0(\partial P_m)+2(\beta_1(S)-1)-m\chi(\partial M).\eqno
(2)$$

{\bf Lemma 5.2}\, $\beta_0(\partial P_m)\le 2\beta_0(S\cap
\partial M)$ for any $m$.

\noindent{\it Proof.}\, The bottom and the top of $P_m$ are
$S^-_1\sqcup S^+_m$, which consists of $2\beta_0(S\cap
\partial M)$ boundary components.
If for some $m$, $\beta_0(\partial P_m)> 2\beta_0(S\cap
\partial M)$, then some component $F$ of $\partial P_m$ does not
meet the top and the bottom of $P_m$. It follows that $F\subset
P_m \subset \widetilde M_S$ provides a component of
$\partial\widetilde M_S$, therefore $p(F)$ is a component of
$\partial M$, where $p: \widetilde M_S\to M$ is the infinite
cyclic covering map. Since the deck transformation group of the
covering $p: \widetilde M_S\to M$ is the infinite cyclic group
which contains no non-trivial finite subgroup, it follows that $p:
F\to p(F)$ is a homeomorphism. Now $S\cap p(F)= \cup_{i=2}^m p(
S_i \cap F)$.

Since $S_i$ separates $P_m$, $S_i$ separates $F$. Since $F$ is
closed, $S_i\cap F$ is homologically trivial in $F$. Hence $p( S_i
\cap F)$ is homologically trivial in $p(F)$, and then $[S\cap
p(F)]=0$, contradicting the assumption in Theorem 5.1 (1).
\hfill\qedsymbol

By using (*) to various homology sequences, we have

\begin{eqnarray*}
\beta_1(Y)&\ge &\beta_1(Y,Y\setminus P_m)-\beta_0(Y\setminus P_m)\qquad\text {by (*)}\\
&= &\beta_1(P_m,\partial P_m)-\beta_0(Y\setminus P_m) \qquad
\text {by excision}\\
&\ge &\beta_1(P_m,\partial P_m)-\beta_0(\partial P_m)
\qquad\;\;\;\;
\text {since $\beta_0(Y\setminus P_m)\le \beta_0(\partial P_m)$}\\
&\ge &\beta_1(P_m)-\beta_1(\partial P_m)-\beta_0(\partial P_m)\qquad\text {by (*)}\\
&\ge &m(\beta_1(X)-\beta_1(S)+\chi(\partial M))+C \qquad \text {by
(1), (2) and Lemma 5.2}
\end{eqnarray*}
where $C=2-\beta_1(S)-6\beta_0(S\cap \partial M)$ is independent
of $m$.

It follows that if $\beta_1(X)>\beta_1(S)-\chi(\partial M)$,
$\beta_1(Y)$ would be arbitrarily large when $m$ gets large. We
reach a contradiction, since $\beta_1(Y)$ should be finite for a
compact manifold $Y$. Theorem 5.1 (1) is proved.

A similar and more direct argument proves Theorem 5.1 (2)
\end{proof}

\vskip 0.3 true cm

{\bf Remark 2.}\;
 Consider the connected sum $M= P\# E(k)$, where $P$
is a homology 3-sphere with $\pi_1(P)\ne 1$ and $k$ is a knot in
$\mathbb S^3$. Let $S\subset M$ be a Seifert surface of $E(k)$,
and $X=M\setminus S$. Then $\beta_1(X)\le \beta_1(S)$ and
$\chi(\partial M)=0$. So the inequality in Theorem 5.1 (1) is not
met.  There is an essential 2-sphere $S^2$ in the connected sum,
and $p^{-1}(S^2)$ is an infinite family of essential 2-spheres  in
$\widetilde M_S$, where $p: \widetilde M_S\to M$ is the infinite
cyclic covering. Then $\widetilde M_S$ can not stay in a compact
3-manifold.

Otherwise suppose $\widetilde M_S\subset Y$ for a compact
3-manifold $Y$. Let $\cup_{i=1}^n S^2_{i}$ be $n$ components in
$p^{-1}(S^2)$ for any given $n$. Then clearly each component of
$Y\setminus \cup_{i=1}^n S_i^2$ contains a copy of the 1-punctured
homology 3-sphere $P^*$ with $\pi_1(P^*)\ne 1$. Since $P^*$ is not
a subset of a punctured 3-sphere, no component of $Y\setminus
\cup_{i=1}^n S_i^2$ is a punctured 3-sphere, which contradicts the
Kneser finiteness theorem [H, Lemma 3.14].

\end{document}